\input amstex

\documentstyle{amsppt}

\rightheadtext{zero-cyles on Severi-Brauer fibrations}

\refstyle{C}

\define\dgr{\roman{deg}\e}

\TagsAsMath

\TagsOnRight

\define\br{\roman{Br}\e}

\define\iv{\roman{Div}\,}

\define\pico{\roman{Pic}^{@,@, 0}}

\define\spec{\roman{Spec}\,}

\define\cb{\overline{C}}

\define\xb{\overline{X}}
\define\xvb{\overline{X}_{\be v}}

\define\chox{\slanted{CH}_{@, 0}@!@!(X)}
\define\choc{\slanted{CH}_{@,@,@, 0}@!@!(C @,@,@,)}
\define\chocb{\slanted{CH}_{@,@,@,0}\ng
\left(\e\cb\e\right)}

\define\choxv{\slanted{CH}_{@,@,@, 0}@!@!(X_{v})}
\define\choxvb{\slanted{CH}_{@,@,@, 0}@!@!(\e\xvb\e)}
\define\chocv{\slanted{CH}_{@,@,@, 0}@!@!(C_{v})}

\define\aox{\slanted{A}_{@,@,@, 0}@!@!(X)}
\define\aoxv{\slanted{A}_{@,@,@, 0}@!@!(\e X_{v}\e)}

\define\choxc{\slanted{CH}_{@,@,@, 0}@!@!(X/C\e)}

\define\choxvcv{\slanted{CH}_{@,@,@, 0}@!@!(X_{v}/C_{v})}

\define\e{@,@,@,@,}
\define\ng{\negthickspace}
\define\be{@!@!@!@!@!}
\define\lbe{@!@!@!}

\define\g{\varGamma}
\define\gv{\varGamma_{v}}

\define\kb{\overline{k}}
\define\kvb{\overline{k_{v}}}

\define\krn{\roman{Ker}\e @,@,}
\define\img{{\roman{Im}}\,}
\define\cok{\roman{Coker}\,}

\define\pic{\roman{Pic}\e\e}
\define\ns{\roman{NS}\!@!@!@!\left(@,@,\xb\e\right)}
\define\nsv{\roman{NS}(@,@,\xvb\e)}
\define\ra{\rightarrow}

\font\tencyr=wncyr10

\font\sevencyr=wncyr7

\font\fivecyr=wncyr5
\newfam\cyrfam
\def\cyr{\fam\cyrfam\tencyr}
\textfont\cyrfam=\tencyr \scriptfont\cyrfam=\sevencyr
\scriptscriptfont\cyrfam=\fivecyr \loadbold \topmatter
\document

\title On the Hasse principle for zero-cycles on
Severi-Brauer fibrations\endtitle

\author Cristian D. Gonz\'alez-Avil\'es\endauthor
\address{Departamento de Matem\'aticas, Universidad
Andr\'es Bello, Sazi\'e 2320, first floor, Santiago,
Chile.}\endaddress

\email{cristiangonzalez\@unab.cl}\endemail

\bigskip

\abstract{Let $k$ be a number field, let $C$ be a smooth,
projective and geometrically integral $k$-curve and let $\pi\:X\ra
C$ be a Severi-Brauer fibration of squarefree index. Various
authors have studied the cokernel of the natural map
$CH_{0}(X/C)\ra\bigoplus_{v}\be CH_{0}(X_{v}/C_{v})$, where
$CH_{0}(X/C)=\krn\!\left[\,\pi_{*}\:CH_{0}(X)\ra
CH_{0}(C)\e\right]$. In this paper we study the kernel of the
above map and find sufficient conditions for it to agree with the
Tate-Shafarevich group of the N\'eron-Severi torus of
$X$.}\endabstract

\footnote"\phantom{.}"{2000 {\it Mathematics Subject
Classification}\;\; 14C15, 14C25, 14K15.}


\endtopmatter

\heading 0. Introduction.\endheading

Let $k$ be a number field, let $C$ be a smooth, projective and
geometrically integral $k$-curve and let $\pi\:X\ra C$ be a
Severi-Brauer fibration of squarefree index. In the case that
$C={\Bbb P}^{1}_{k}$ and $X\ra C$ is a conic bundle,
J.-L.Colliot-Th\'el\`ene and J.-J.Sansuc conjectured in 1981[2]
the existence of an exact sequence of torsion groups
$$
0\ra{\cyr X}^{\e 1}\be(T\e)\ra\aox\ra\bigoplus_{v}\aoxv\ra
\roman{Hom}\!\left(\e H^{@,@,1}\ng\left(k,\ns\right),\Bbb Q/\Bbb
Z\e\right),\tag 1
$$
where $T$ is the N\'eron-Severi torus of $X$ and $\aox$ is the
group of rational equivalence classes of 0-cycles of degree zero
on $X$. This conjecture was proved in 1988 by P.Salberger, in his
remarkable paper [10]. Then, in 1999 [3], J.-L.Colliot-Th\'el\`ene
partially extended Salberger's proof to conic bundles $X\ra C$
over a curve $C$ of arbitrary genus. Under certain extra
assumptions which turn out to be unnecessary, this author obtained
an exact sequence of torsion groups
$$
\choxc\ra\bigoplus_{v}\choxvcv\ra \roman{Hom}(\e\br X/\br C,\Bbb
Q/\Bbb Z\e),
$$
where $\choxc=\krn\!\left[\,\pi_{*}\:\chox\ra\choc\e\right]$ is
the relative 0-th Chow group of the fibration $X\ra C$. On the
other hand, in his 1990 ICM talk [1], S.Bloch put forth a general
conjecture which, in the case of conic bundles over a curve $C$,
asserts the equality of the groups
$$
{\cyr X}\choxc=\krn\ng\left[\choxc\ra\prod_{v}\choxvcv\right]
$$
and ${\cyr X}^{\e 1}\be(T\e)$. This direct generalization of part
of Salberger's result was shown to be incorrect by V.Suresh in his
1997 paper [11]. This author produced an example of a conic bundle
over an elliptic curve over $\Bbb Q\,$\footnote{Suresh's base
curve has multiplicative reduction at 3 and 7, so neither of the
hypotheses of Corollary 0.2 is satisfied in his example.} for
which ${\cyr X}^{\e 1}\be(T\e)=0$ but ${\cyr X}\choxc\neq 0$.

In this paper we study the group ${\cyr X}\choxc$ for an arbitrary
Severi-Brauer fibration $X\ra C$ of squarefree index. Our main
result is the following.

\proclaim{Theorem 0.1} Let $k$ be a number field, let $X\ra C$ be
a Severi-Brauer $k$-fibration of squarefree index and let
$\Phi\:\choxc\ra H^{@,@,@,1}(k,T)$ be the characteristic map of
$X\ra C$, where $T$ is the N\'eron-Severi torus of $X$. Then there
exists a natural exact sequence
$$
0\ra{\cyr X}(\e\krn\Phi)\ra{\cyr X}\choxc\ra{\cyr X}^{\e
1}\be(T)\ra 0.
$$
\endproclaim

In Section 2 we use a well-known Hasse-principle theorem due to
K.Kato to obtain the following corollary of Theorem 0.1.

\proclaim{Corollary 0.2} Let $k$ be a number field and let $X\ra
C$ be a Severi-Brauer $k$-fibration of squarefree index $n$.
Assume that the following conditions hold:

\roster

\item"i)" either $k$ is totally imaginary or $n$ is odd, and

\item"ii)" $C$ acquires good reduction at all primes over an
extension of $k$ of degree prime to $n$.

\endroster

Then there exists a canonical isomorphism
$$
{\cyr X}\choxc={\cyr X}^{\e 1}\be(T\e),
$$
where $T$ is the N\'eron-Severi torus of $X$.
\endproclaim

The key ingredient of the proof of Theorem 0.1 is an approximation
lemma of E.Frossard [7, Lemma 3.3]. This fundamental result is
used in the proofs of Lemmas 2.1 and 2.3 below, which form the
technical center of the paper. These lemmas, together with some
elementary homological algebra, yield Theorem 0.1 above.

Combining Corollary 0.2 and [7, Theorem 0.4], we obtain the
following result.

\proclaim{Corollary 0.3} Let the hypotheses and notations be as in
the statement of Corollary 0.2. Then there exists a natural exact
sequence of torsion groups
$$
0\ra{\cyr X}^{\e 1}\be(T\e)\ra\choxc\ra\bigoplus_{v}\choxvcv\ra
\roman{Hom}(\e\br X/\br C,\Bbb Q/\Bbb Z\e).
$$
\endproclaim

We also obtain the following theorem, which generalizes
Salberger's exact sequence (1)\footnote{This result must be
familiar to all specialists in this area. See Remark 4.9 of [5],
where the existence of such an exact sequence is suggested (but
not formally stated).}

\proclaim{Theorem 0.4} Let $k$ be a number field and let $X\ra
\Bbb P^{1}_{k}$ be a Severi-Brauer $k$-fibration of squarefree
index. Assume that $X\be(k_{v})\neq\emptyset$ for every real prime
$v$ of $k$. Then there exists a natural exact sequence
$$
0\ra{\cyr X}^{\e 1}\be(T\e)\ra A_{0}(X)\ra \bigoplus_{v}
A_{0}(X_{v})\ra\roman{Hom}\!\left(\e
H^{@,@,1}\ng\left(k,\ns\right),\Bbb Q/\Bbb Z\e\right),
$$
where $T$ is the N\'eron-Severi torus of $X$.\qed
\endproclaim

 \heading Acknowledgement\endheading

I thank E.Frossard and J.van Hamel for sending me copies of their
papers [6,7,12]. As mentioned above, the work of Frossard is crucial
to this paper, and I wish to express my gratitude to her here for
having laid the foundations on which this paper rests. Finally, I
thank the referees for several helpful suggestions.

\heading 1. Preliminaries.\endheading

Let $k$ be a field, fix a separable closure $\kb$ of $k$ and let
$\g=\roman{Gal}\!\lbe\left(\e\e\kb/k\right)$. If $M$ is a discrete
and continuous $\g$-module, the $i$-th $\g$-cohomology group of $M$
will be denoted by $H^{i}(k,M)$. Now assume that $k$ is a number
field. For each prime $v$ of $k$, we will write $k_{v}$ for the
completion of $k$ at $v$. We define
$$
{\cyr X}^{\e i}\be(M)=\krn\ng\left[H^{i}(k,M)\ra\prod_{\text{all}\,
v} H^{i}(k_{v},M)\right],
$$
where the map involved is the natural localization map.

Now let $X$ be a smooth, projective and geometrically integral
$k$-variety. We will write $\xb$ for the $\kb$-variety
$X\otimes_{k}\kb$, $k(X)$ (resp. $\kb(X)$) for the field of
rational functions on $X$ (resp. $\xb$) and $\ns$ for the
N\'eron-Severi group of $\xb$. In addition, the group of 0-cycles
on $X$ will be denoted by $Z_{0}(X)$ and $\chox$ will denote its
quotient by rational equivalence. Further, we will write
$A_{0}(X)$ for the kernel of the degree map $\dgr\:\chox\ra \Bbb
Z$ and, for each prime $v$ of $k$, $X_{v}$ will denote the
$k_{v}$-variety $X\otimes_{k}k_{v}$.

Let $C$ be a smooth, projective and geometrically integral
$k$-curve and let $X\ra C$ be a Severi-Brauer $k$-fibration, i.e.,
$X$ is a smooth, projective and geometrically integral $k$-variety
equipped with a proper and dominant $k$-morphism $\pi\: X\ra C$
whose generic fiber $X_{\eta}=X\times_{C}\spec k(C)$ is a
Severi-Brauer $k(C)$-variety.

Given a central simple $k(C)$-algebra $A$, there exists an Artin
model $X$ of $A$ over $C$ [6]. Thus $X$ is a smooth, projective
and geometrically integral $k$-variety equipped with a proper and
dominant $k$-morphism $\pi\: X\ra C$ whose generic fiber
$X_{\eta}=X\times_{C}\spec k(C)$ is the Severi-Brauer
$k(C)$-variety associated to $A$. We will write $X_{\kb(\eta)}$
for the generic fiber of $\overline{\pi}\:\xb\ra\cb$. Note that,
since $\kb(C)$ is a $C_{1}-$field, $X_{\kb(\eta)}$ is a projective
space and therefore $\pic X_{\kb(\eta)}=\Bbb Z$. The {\it index}
\, of such a fibration is defined to be the index of $A$, i.e.,
the minimum degree of all extensions of $k(C)$ which split $A$ or,
equivalently, the least positive degree of a 0-cycle on
$X_{\eta}$. We are interested in the group
$$
{\cyr X}\chox=\krn\ng\left[\chox\ra \prod_{v}\choxv\right].\tag 2
$$
Set
$$
\choxc=\krn\!\left[\,\pi_{*}\:\chox\ra\choc\e\right].
$$
Since the map $\choc\ra\prod\chocv$ may be identified with the map
$\pic C\ra\prod\pic C_{v}$, which is injective (see, e.g., [4],
Proposition 1.1, p.3), we have
$$
{\cyr X}\chox={\cyr X}\choxc\,,
$$
where
$$
{\cyr X}\choxc=\krn\ng\left[\choxc\ra
\prod_{v}\choxvcv\right].\tag 3
$$
Thus the study of (2) is reduced to the study of (3). Assume now
that $X\ra C$ has squarefree index. Then (2) (and hence also (3))
is a finite group [6, Th\'eor\`eme 4.7]. Further, $\choxvcv$ is
zero for all but finitely many primes $v$ of $k$ [op.cit.,
Th\'eor\`eme 4.8]. We will write $S$ for the set of primes $v$ of
$k$ such that $\choxvcv\neq 0$.

We now observe that, in studying the group (3), we may restrict
our attention to the chosen model $X\ra C$ of $A$, for if
$X^{\prime}\ra C$ is an arbitrary Severi-Brauer $k$-fibration
whose generic fiber corresponds to $A$, then $X^{\prime}\ra C$ is
birationally equivalent to $X$ over $C$ (as follows from [EGA IV,
8.8.2.5]), and the Chow group of dimension zero is a birational
invariant of smooth projective varieties (see [8, 16.1.11, p.312],
and note that the proof given there holds over any field).

\medpagebreak

We now write $C_{0}$ for the set of closed points of $C$. For each
$P\in C_{0}$, with residue field $k(P)$, we let
$X_{P}=X\times_{C}\spec k(P)$ be the fiber of $\pi$ over $P$. We
will write $\Cal P=\{P_{i}\}_{i\in I}$ for the finite set of
closed points of $C$ for which the fiber $X_{P}$ is singular
(i.e., not geometrically integral). Further, for each $i\in I$, we
will write $k_{i}$ for $k(P_{i})$. For each $P\in C_{0}$, there
exists a natural residue homomorphism
$$
\roman{res}_{\e P}\: H^{@,@,2}\lbe(k(C),\mu_{n})\ra H^{1}(k(P),\Bbb
Z/n\e\Bbb
Z)=\roman{Hom}\!\left(\roman{Gal}\!\left(\e\kb/k(P)\right),\Bbb
Z/n\e\Bbb Z\right),
$$
where $n$ is the index of $A$ (see [6, pp.205-206]). The algebra
$A$ defines an element of $H^{@,@,2}\lbe(k(C),\mu_{n})$ and, for
each $i\in I$, we will write $L_{i}$ for the cyclic extension of
$k_{i}$ determined by the residue $\roman{res}_{\e P_{\be
i}}\be(A)$, i.e., $L_{i}$ is the fixed field of the kernel of
$\roman{res}_{\e P_{\be i}}\be(A)$. Let $E=\prod_{i\in I}k_{i}$,
$F=\prod_{i\in I}L_{i}$, $N_{i}=N_{L_{i}/k_{\e i}}$ and
$N=N_{F/E}$. The quotient of $E^{\e *}/NF^{\e *}=\prod_{i\in
I}k_{i}^{*}/N_{i}L_{i}^{*}$ by the image of the diagonal map
$k^{*}\ra E^{\e *}$ will be denoted by $E^{\e *}/k^{*}NF^{\e *}$.
Clearly, $E^{\e *}/k^{*}NF^{\e *}$ is annihilated by $n$. Now, for
each prime $v$ of $k$ (and each $i\in I$), we will write $k_{i,v}$
(resp. $L_{i,v}$, $E_{v}$, $F_{v}$) for $k_{\e i}\otimes_{k}k_{v}$
(resp. $L_{i}\otimes_{k}k_{v}$, $E\otimes_{k}k_{v}$,
$F\otimes_{k}k_{v}$). Further, $A_{v}$ will denote
$A\otimes_{k(C)}k_{v}\be(C)$. There exists a natural exact
sequence
$$
0\ra B_{1}\ra B_{2}\be\oplus\Bbb Z\ra\ns\ra\Bbb Z\ra 0,\tag 4
$$
where the map $\ns\ra\Bbb Z$ is given by restricting divisors to
the generic fiber, $B_{1}$ is the free group on the closed points
$y$ of $\cb$ for which the fiber of $\overline{\pi}$ over $y$,
$\xb_{\be y}=\xb\times_{\cb} \e\spec\kb(y)$, is singular,
$B_{2}\subset\iv\xb$ is the free abelian group on the irreducible
components of the singular fibers of $\overline{\pi}$ and the map
$B_{1}\ra B_{2}\be\oplus\Bbb Z$ is induced by $y\mapsto
\left(\e\xb_{\be y},-1\right)$ (see [6, p.201] and note that
$\overline{\pi}^{\e *}\pico\e\cb=\pico\e\xb$ [7]). The groups
$B_{1}$ and $B_{2}$ are isomorphic to $\bigoplus_{i\in I}\Bbb
Z[\roman{Gal}(k_{i}/k)]$ and $\bigoplus_{i\in I}\Bbb
Z[\roman{Gal}(L_{i}/k)]$, respectively. It is not difficult to
check that $(B_{2}\be\oplus\Bbb Z)/B_{1}$ is torsion-free. Let
$T_{0}$ denote the $k$-torus corresponding to $(B_{2}\be\oplus\Bbb
Z)/B_{1}$, and let $T$ be the N\'eron-Severi torus of $\xb$. Then
there exists a natural isomorphism
$$
\delta\: E^{\e *}/k^{*}NF^{\e *}\ra H^{@,@,@,1}\be(k,T_{0}).\tag 5
$$
See [6, Lemma 3.5]. Further, (4) induces a natural exact sequence
$$
0\ra H^{@,@,@,1}\be(k,T\e)\ra H^{@,@,@,1}\be(k,T_{0})\ra\br k.
$$
A similar exact sequence exists over $k_{v}$ for every prime $v$ of
$k$, and we have a natural exact commutative diagram
$$\minCDarrowwidth{.8cm}
\CD 0 @>>>H^{@,@,@,1}\be(k,T\e)  @>>>
H^{@,@,@,1}\be(k,T_{0})@>>>\br k\\
      @.   @VVV     @VVV     @VVV\\
      0 @>>>\underset{v}\to\bigoplus\,H^{@,@,@,1}\be(k_{v},T\e)
      @>>>\underset{v}\to\bigoplus\,
H^{@,@,@,1}\be(k_{v},T_{0}) @>>>\underset{v}\to\bigoplus\,\br k_{v},
\endCD
$$
where the direct sums extend over all primes $v$ of $k$. Now the
commutativity of the above diagram and the classical
Brauer-Hasse-Noether theorem yield

\proclaim{Lemma 1.1} There exists a canonical isomorphism ${\cyr
X}^{\e 1}\be(@,@,T@,@,)={\cyr X}^{\e 1}\be(@,@,T_{0})$.
\endproclaim

\medpagebreak

Recall that ${\cyr X}^{@,@,1}\be(\e T\e)$ is a birational
invariant of $X$.

\medpagebreak

We will need the following definition.

\proclaim{Definition}{\rom{The group of {\it divisorial norms},
denoted $k(C)_{\roman{dn}}^{*}$, is the subgroup of $k(C)^{*}$
consisting of functions which, at every point $P$ of $C$, can be
written as the product of a unit at $P$ by an element of the
reduced norm group $\roman{Nrd}\e A^{*}$. Equivalently}}
$$
k(C)_{\roman{dn}}^{*}=\{ f\in
k(C)^{*}\:\roman{div}_{C}f\in\pi_{*}(Z_{0}(X))\}.
$$
\endproclaim
There exists a canonical isomorphism
$$
\Psi\:\choxc\ra k(C)_{\roman{dn}}^{*}\big/k^{*}\e\roman{Nrd}\e
A^{*}\tag 6
$$
which maps a cycle $z\in\choxc$ to the function $f$ (defined up to
multiplication by a constant) such that
$\roman{div}_{C}f=\pi_{*}\be(z)$. See [6, p.98]. Now, by
evaluating functions at the points of $\Cal P$, one obtains a
natural ``specialization map"
$$
\gamma\:k(C)_{\roman{dn}}^{*}/k^{*}\e\roman{Nrd}A^{*}\ra
E^{*}/k^{*}NF^{*},\tag 7
$$
and E.Frossard [6] has defined a ``characteristic" map
$$
\Phi_{0}\:\choxc\ra H^{@,@,@,1}(k,T_{0}).\tag 8
$$
These functions (5)-(8) fit into a natural commutative diagram
$$
\CD \choxc  @>\Psi>\sim> k(C)_{\roman{dn}}^{*}/k^{*}\e
\roman{Nrd}A^{*}\\
@V\Phi_{0}VV     @V\gamma VV   \\
H^{@,@,@,1}(k,T_{0})@<\delta<\sim< E^{*}/k^{*}NF^{*}.
\endCD\tag 9
$$
See [6, Proposition 3.8].

\proclaim{Lemma 1.2} Let $k$ be a number field and let $X\ra C$ be
an Artin model over $C$ of a central simple $k(C)$-algebra $A$ of
squarefree index. Then the map (6) induces an isomorphism
$\krn\Phi_{0}=\krn\gamma$.
\endproclaim
{\it Proof.} See Corollary 3.9 and Proposition 4.1 of [6].\qed

\medpagebreak

\heading 2. The proofs.\endheading

In this Section we establish the main results of the paper. We
keep the notations and hypotheses introduced in the previous
Section.

\proclaim{Lemma 2.1} Let $A$ be a central simple $k(C)$-algebra of
squarefree index, let $X\ra C$ be an Artin model of $A$ over $C$
and let $\Phi_{0}$ be the characteristic map (8). Then the natural
map
$$
\krn\e\Phi_{0}\ra\prod_{v}\krn\e\Phi_{0,v}
$$
is surjective.
\endproclaim
{\it Proof.} By Lemma 1.2, it suffices to show that the natural
map $\krn\e\gamma\ra\prod_{v}\krn\e\gamma_{v}$ is surjective.
Using a transfer argument analogous to that used in [7, pp.94-95],
we may reduce the proof to the case where the fibration $X\ra C$
satisfies hypothesis 2.1 of [op.cit.]\footnote{With the notations
of [op.cit.], the cokernel of $\krn\e\gamma\ra\prod_{\,
v}\be\krn\e\gamma_{v}$ is annihilated by multiplication by the
index of $X\ra C$ as well as by $[M_{i}\:k]$ for each $i$ if the
assertion of the lemma is true for the fibration
$X_{M_{i}^{\prime}}\ra C_{M_{i}^{\prime}}$ of [op.cit.], p.94,
line -1. This immediately yields the surjectivity of
$\krn\e\gamma\ra\prod_{v}\krn\e\gamma_{v}$.}. Now recall the set
$S$ of primes $v$ of $k$ such that $\choxvcv\neq 0$. Clearly if
$v\notin S$ then $\krn\gamma_{v}=\krn\Phi_{0,v}=0$, whence
$\prod_{v}\krn\e\gamma_{v}=\prod_{v\in S}\krn\e\gamma_{v}$. Let
$(f_{v})_{v\in S}\in\prod_{v\in S}\krn\e\gamma_{v}$. For each
$v\in S$, lift $f_{v}$ to $\tilde{f}_{v}\in
k_{v}(C)^{*}_{\roman{dn}}/k_{v}^{*}$ and choose $z_{v}\in
Z_{0}(X_{v})$ such that
$(\pi_{v})_{*}(z_{v})=\roman{div}_{C_{v}}\tilde{f}_{v}$, where
$\pi_{v}=\pi\otimes_{k} k_{v}$. Then Lemma 3.3 of [7] (which may
be applied since $X\ra C$ was assumed to satisfy hypothesis 2.1 of
[op.cit.]), with $y=0$, $\eta=1\in E^{\e *}/k^{*}NF^{\e *}$,
$f_{z_{v},y}=\tilde{f}_{v}$ for $v\in S$ and $f_{z_{v},y}=1$ for
$v\notin S$, asserts the existence of an element $g\in\krn\gamma$
which maps to $(f_{v})_{v\in S}$, as desired (see Remark 2.2
below).\qed

\medpagebreak

{\bf Remark 2.2.} Regarding the proof of the lemma, the following
comment is in order. The element $g\in k(C)^{*}\be/k^{*}$ whose
existence is asserted by Lemma 3.3 of [7] may be chosen so that it
specializes to the given $\eta\in E^{\e *}/k^{*}NF^{\e *}$. This
follows from the fact that the function $h\in k(C)^{*}$ which
appears in the proof of that lemma [op.cit.], p.94, line 16, and
whose existence is asserted by Proposition 1.4 of [op.cit.], may
be chosen so that it specializes to $\hat{\eta}^{\prime}$, as
follows from the {\it proof} of Proposition 1.4 of [op.cit.]. See
[op.cit.], p.90, line 9.

\medpagebreak

\proclaim{Lemma 2.3} Let $A$ be a central simple $k(C)$-algebra of
squarefree index, let $X\ra C$ be an Artin model of $A$ over $C$
and let $\Phi_{0}$ be the characteristic map (8). Then the natural
map
$$
\cok\Phi_{0}\ra\prod_{v}\cok\Phi_{0,v}
$$
is injective.
\endproclaim
{\it Proof.} As in the proof of Lemma 2.1, a transfer argument
analogous to that used in [7, pp.94-95] enables us to reduce the
proof to the case where $X\ra C$ satisfies hypothesis 2.1 of
[7]\footnote{The essential fact again is that, for any (Galois)
extension $M/k$, the kernel of the map
$\cok\Phi_{0}\ra\cok\Phi_{\e 0,M}$ induced by the restriction map
is annihilated both by multiplication by $[M\: k]$ and by
multiplication by the index of $X\ra C$.}. Let $\alpha\in
H^{@,@,1}\be(k,T_{0})$ be such that $\roman{res}_{\e
v}\be(\alpha)\in\img\Phi_{0,v}$ for every prime $v$ of $k$, where
$\roman{res}_{v}\:H^{@,@,1}\be(k,T_{0})\ra
H^{@,@,1}\be(k_{v},T_{0})$ is the natural restriction map. For
each $v$, choose $z_{v}\in\krn[\e (\be\pi_{v}\be)_{*}\:
Z_{0}(X_{v})\ra\chocv]$ such that $\roman{res}_{\e
v}\be(\alpha)=\Phi_{0,v}\!\left([z_{v}]\right)$, where $[z_{v}]$
denotes the class of $z_{v}$ in $\choxvcv$. Further, let $f_{v}\in
k_{v}(C)^{*}\be/k_{v}^{*}$ be such that
$(\be\pi_{v}\be)_{*}(z_{v})=\roman{div}_{C_{v}}\be(f_{v})$. Now
let $\eta=\delta^{-1}(\alpha)\in E^{\e *}/k^{*}NF^{\e *}$, where
$\delta$ is the map (5). For each $v$, the commutativity of
diagram (9) (over $k_{v}$) shows that $f_{v}$ specializes to
$\delta_{v}^{-1}(\Phi_{0,v}([z_{v}]))$, where $\delta_{v}$ is the
map (5) for $k_{v}$, i.e., to the image of $\eta$ in $E_{v}^{\e
*}/k_{v}^{*}\e NF_{v}^{\e *}$. Now Lemma 3.3 of [7], with $y=0$
and $\eta$ as above, asserts the existence of elements $g\in
k(C)^{*}_{\roman{dn}}\big/k^{*}\roman{Nrd}\e A^{*}$ and $z\in
Z_{0}(X)$ such that $g$ specializes to $\eta$ (see Remark 2.2
above) and $\Psi([z])=g$, where $\Psi$ is the map (6). It follows
that $\delta^{-1}(\alpha)=\delta^{-1}(\Phi_{0}([z]))$ (see diagram
(9)), whence $\alpha=\Phi_{0}([z])$. This completes the proof.
\qed

\medpagebreak

{\bf Remark 2.4.} In the case of conic bundles, Lemmas 2.1 and 2.3
are immediate consequences of Theorem 1.3 of [3] (to obtain Lemma
2.1 from loc.cit., set $\alpha=0$ there). This observation was the
starting point of this paper.

\medpagebreak

We can now state our main result.

\proclaim{Theorem 2.5} Let $k$ be a number field, let $X\ra C$ be
a Severi-Brauer $k$-fibration of squarefree index and let
$\Phi_{0}$ be the characteristic map (8). Then there exists a
natural exact sequence
$$
0\ra{\cyr X}(\e\krn\Phi_{0})\ra{\cyr X}\choxc\ra{\cyr X}^{\e
1}\be(T)\ra 0,
$$
where $T$ is the N\'eron-Severi torus of $X$.
\endproclaim
{\it Proof}. Assume first that $X\ra C$ is the Artin model of a
central simple $k(C)$-algebra $A$ of squarefree index. There
exists a natural exact commutative diagram
$$\minCDarrowwidth{.8cm}
\CD 0 @>>>\krn\Phi_{0}  @>>>\choxc @>>>\img\Phi_{0} @>>>0\\
      @.   @VVV     @VVV     @VVV    @.\\
      0 @>>>\underset{v}\to\prod\,\krn\Phi_{0,v}
      @>>>\underset{v}\to\prod\,
\choxvcv @>>>\underset{v}\to\prod\,\img\Phi_{0,v} @>>>0.
\endCD
$$
Applying the snake lemma to the above diagram and using Lemma 2.1,
we obtain an exact sequence $0\ra{\cyr X}(\e\krn\Phi_{0})\ra{\cyr
X}\choxc\ra {\cyr X}(\e\img\Phi_{0})\ra 0$. On the other hand, the
commutativity of the natural diagram
$$\minCDarrowwidth{.8cm}
\CD 0 @>>>\img\Phi_{0}  @>>> H^{@,@,1}\be(k,T_{0})@>>>\cok\Phi_{0}
@>>>0\\
      @.   @VVV     @VVV     @VVV    @.\\
      0 @>>>\underset{v}\to\prod\,\img\Phi_{0,v}
      @>>>\underset{v}\to\prod\,
      H^{@,@,1}\be(k_{v},T_{0})@>>>\underset{v}\to\prod\,
      \cok\Phi_{0,v}@>>>0
\endCD
$$
together with Lemmas 2.3 and 1.1 yield a natural isomorphism
${\cyr X}(\e\img\Phi_{0})={\cyr X}^{\e 1}\be(T)$. This completes
the proof when $X\ra C$ is an Artin model as above. The general
case may be deduced from the preceding one by using the birational
invariance of the groups ${\cyr X}\choxc$ and ${\cyr X}^{\e
1}\be(T)$.\qed

\medpagebreak

By an example of V.Suresh [11], the group ${\cyr
X}(\e\krn\Phi_{0})$ appearing in the statement of the theorem need
not vanish. We will now study this group and find conditions under
which it vanishes (see Theorem 2.8 below).

\medpagebreak

There exists a natural exact commutative diagram
$$\minCDarrowwidth{.6cm}
\CD 0 @>>>k^{*}/\e k^{*}\cap\roman{Nrd}\e A^{*} @>>>
k(C)_{\roman{dn}}^{*}\big/\e\roman{Nrd}\e A^{*}@>>>k(C)_{
\roman{dn}}^{*}\big/k^{*}\e\roman{Nrd}\e A^{*}@>>>0\\
      @.   @V\widetilde{\gamma}_{\e 0}VV     @V\widetilde{\gamma} VV   @V\gamma VV   @.\\
      0 @>>>k^{*}/\e k^{*}\cap NF^{\e *}@>>>
E^{\e *}\be/NF^{\e *}@>>>E^{\e *}/k^{*}NF^{\e *}@>>>0,
\endCD \tag 10
$$
where $\widetilde{\gamma}_{0}$ is induced by the identity on
$k^{*}$ and $\widetilde{\gamma}$ and $\gamma$ and induced by
specialization. Since $\widetilde{\gamma}_{0}$ is surjective, the
snake lemma applied to (10) yields a natural exact sequence
$$
0\ra\krn\widetilde{\gamma}_{\e 0}\ra\krn\widetilde{\gamma}\ra
\krn\gamma\ra 0.\tag 11
$$
Note that $\krn\widetilde{\gamma}_{\e 0}$ is canonically
isomorphic to $k^{*}\cap(\cap_{\e i\in I}\e
N_{i}L_{i}^{*})\e\big/\e k^{*}\cap\roman{Nrd}\e A^{*}$.

\medpagebreak

Now there exist sequences analogous to (11) over $k_{v}$ for every
prime $v$ of $k$, and we have a natural exact commutative diagram
$$\minCDarrowwidth{.6cm}
\CD 0 @>>>\krn\widetilde{\gamma}_{\e 0}@>>>
\krn\widetilde{\gamma}@>>>\krn\gamma@>>>0\\
      @.   @VVV     @VVV   @VVV   @.\\
      0 @>>>\underset{v}\to\prod\,\krn\widetilde{\gamma}_{\e 0,v}
      @>>>\underset{v}\to\prod\,\krn\widetilde{\gamma}_{v}
      @>>>\underset{v}\to\prod\,\krn\gamma_{v} @>>>0.
\endCD\tag 12
$$

\medpagebreak

We now have\footnote{One of the referees informed us that this
result, as well as Proposition 2.7 and Theorem 2.9 below, have
been known to the specialists in this area for quite some time.
Unfortunately, we have been unable to find appropriate references
for them and therefore provide the corresponding proofs.}

\proclaim{Lemma 2.6} The canonical map
$$
\krn\e\widetilde{\gamma}\ra\prod_{v}\krn\e\widetilde{\gamma}_{v}
$$
is injective.
\endproclaim
{\it Proof}. P.Salberger (see [6, Proposition 4.4]) has shown
that, over any field $k$ of characteristic 0, there exists a
canonical injection $\krn\e\widetilde{\gamma}\hookrightarrow
H_{\roman{nr}}^{@,@,3}\ng\left(\e k(C)/C,\mu_{n}^{\otimes
2}\e\right)$, where
$$
H_{\roman{nr}}^{@,@,3}\ng\left(\e k(C)/C,\mu_{n}^{\otimes
2}\e\right)=\krn\ng \left[H^{@,@,3}\ng\left(\e
k(C),\mu_{n}^{\otimes 2}\e\right)\ra \bigoplus_{P\in C_{0}}
H^{@,@,2}\ng\left(\e k(P),\mu_{n}\e\right)\right].
$$
Here, for each $P\in C_{0}$, $H^{@,@,3}\!\be\left(\e
k(C),\mu_{n}^{\otimes\e 2}\e\right)\ra H^{@,@,2}\!\left(\e
k(P),\mu_{n}\e\right)$ is the natural residue map. On the other
hand, a well-known theorem of K.Kato [9, \S 4], asserts that the
canonical map
$$
H_{\roman{nr}}^{@,@,3}\ng\left(\e k(C)/C,\mu_{n}^{\otimes
2}\e\right)\ra\bigoplus_{v}H_{\roman{nr}}^{@,@,3}\ng\left(\e
k_{v}(C)/C_{v},\mu_{n}^{\otimes 2}\e\right)
$$
is an isomorphism, whence the lemma follows.\qed

\medpagebreak

We now apply the snake lemma to diagram (12) and use the preceding
lemma together with Lemma $1.2\e$ \footnote{As well as a
``birational invariance" argument as in the proof of Theorem 2.4.}
to obtain the following result

\proclaim{Proposition 2.7} Let $k$ be a number field, let $X\ra C$
be a Severi-Brauer $k$-fibration of squarefree index and let
$\Phi_{0}$ be the characteristic map (8). Then there exists a
canonical injection
$$
{\cyr X}(\e\krn\Phi_{0})\hookrightarrow\cok\!\Bigg[
\krn\widetilde{\gamma}_{0}\ra\prod_{v}\krn\widetilde{\gamma}_{0,v}\Bigg],$$
where $\widetilde{\gamma}_{0}$ is the map in diagram (10). \qed
\endproclaim

\medpagebreak

\proclaim{Theorem 2.8} Let $k$ be a number field, let $X\ra C$ be
a Severi-Brauer $k$-fibration of squarefree index with associated
central simple $k(C)$-algebra $A$, and let $\Phi_{0}$ be the
characteristic map (8). Further, let $S^{\e\prime}$ be the set of
primes of $k$ consisting of the real primes and the primes of bad
reduction for $C$. Assume that
$$
k_{v}^{*}\cap(\cap_{\e i\in I}\e
N_{i,v}L_{i,v}^{*})=k_{v}^{*}\cap\roman{Nrd}\,
A_{v}^{*}\quad\text{for all $v\in S'$}.
$$
Then the natural map $\,\krn\Phi_{0}\ra\prod_{v}\krn\Phi_{0,v}$ is
injective. Consequently, there exists a canonical isomorphism
$$
{\cyr X}\choxc={\cyr X}^{\e 1}\be(T).
$$
\endproclaim
{\it Proof.} Let $S^{\e\prime}$ be as in the statement. Then
$H_{\roman{nr}}^{@,@,3}\!\left(\e k_{v}(C)/C_{v},\mu_{n}^{\otimes
2}\e\right)=0$ for $v\notin S^{\e\prime}$. See [9, Corollary 2.9]
and note that $k_{v}(C)$ has cohomological dimension 1 if $v$ is
complex. Now the proof of Lemma 2.6 shows that
$\krn\widetilde{\gamma}_{v}$, and therefore also
$\krn\widetilde{\gamma}_{\e 0,v}$ (see (11)), vanishes if $v\notin
S^{\e\prime}$. On the other hand, the hypothesis implies that
$\krn\widetilde{\gamma}_{\e 0,v}=0$ for $v\in S^{\e\prime}$. We
conclude that $\prod_{v}\krn\widetilde{\gamma}_{0,v}=0$, and the
theorem now follows from Theorem 2.4 and Proposition 2.7.\qed

\medpagebreak

The next result gives conditions under which the full group
$\krn\Phi_{0}$ (and not just its subgroup ${\cyr
X}(\e\krn\Phi_{0})$) vanishes.

\proclaim{Theorem 2.9} Let $k$ be a number field, let $A$ be a
central simple $k(C)$-algebra of squarefree index $n$ and let
$X\ra C$ be an Artin model of $A$ over $C$. Assume that the
following conditions hold:

\roster

\item"i)" either $k$ is totally imaginary or $n$ is odd, and

\item"ii)" $C$ acquires good reduction at all primes over an
extension of $k$ of degree prime to $n$.

\endroster
Then the characteristic map $\,\Phi_{0}\:\choxc\ra
H^{@,@,@,1}(k,T_{0})$ is injective.
\endproclaim
{\it Proof.} Assume first that $C$ has good reduction at all
primes of $k$. Then $H_{\roman{nr}}^{@,@,3}\!\left(\e
k_{v}(C)/C_{v},\mu_{n}^{\otimes 2}\e\right)=0$ for all
non-archimedean $v$. See [9, Corollary 2.9]. We conclude (cf.
proof of Lemma 2.6) that $\krn\widetilde{\gamma}_{v}=0$ for all
such $v$. The same conclusion holds if $v$ is a complex prime (cf.
proof of Theorem 2.8). If $k_{v}\simeq\Bbb R$, then
$H_{\roman{nr}}^{@,@,3}\!\left(\e k_{v}(C)/C_{v},\mu_{n}^{\otimes
2}\e\right)$ is annihilated by multiplication by $2$ as well as by
multiplication by $n$, so (i) shows that this group is zero as
well. We conclude that $\krn\widetilde{\gamma}_{v}=0$ for all $v$,
whence $\krn\widetilde{\gamma}=0$ by Lemma 2.6. Now (11) shows
that $\krn\gamma=0$, and Lemma 1.2 completes the proof in this
case. Now let $M$ be an extension of $k$ of degree prime to $n$
over which $C$ acquires good reduction at all primes. The first
part of the proof shows that $H_{\roman{nr}}^{@,@,3}\!\left(\e
M_{w}(C)/C_{w},\mu_{n}^{\otimes 2}\e\right)=0$ for all primes $w$
of $M$, where $C_{w}=C\otimes_{k}M_{w}$. On the other hand, if $v$
is a prime of $k$ and $w$ is a prime of $M$ lying above $v$, the
kernel of the restriction map
$$
H_{\roman{nr}}^{@,@,3}\!\left(\e k_{v}(C)/C_{v},\mu_{n}^{\otimes
2}\e\right)\ra H_{\roman{nr}}^{@,@,3}\!\left(\e
M_{w}(C)/C_{w},\mu_{n}^{\otimes 2}\e\right)
$$
is annihilated by multiplication by $[M\:k\e]$ as well as by
multiplication by $n$. We conclude that
$H_{\roman{nr}}^{@,@,3}\!\left(\e k_{v}(C)/C_{v},\mu_{n}^{\otimes
2}\e\right)=0$ for all $v$, and the proof goes through as
before.\qed

\medpagebreak

Now let $\br X/\e\br C$ denote $\br X/\pi^{*}\br C$. Combining
Theorem 0.4 of [7]\footnote{See also the comments following (0.b)
on p.83 of [7].} (which yields the last three terms of the exact
sequence appearing in the next Corollary) and Theorem 2.4, Theorem
2.9 and a ``birational invariance" argument, we obtain

\proclaim{Corollary 2.10} Let $k$ be a number field and let $X\ra
C$ be a Severi-Brauer $k$-fibration of squarefree index $n$.
Assume that the following conditions hold:

\roster

\item"i)" either $k$ is totally imaginary or $n$ is odd, and

\item"ii)" $C$ acquires good reduction at all primes over an
extension of $k$ of degree prime to $n$.

\endroster

Then there exists a natural exact sequence
$$
0\ra{\cyr X}^{\e 1}\be(T\e)\ra\choxc\ra\bigoplus_{v}\choxvcv\ra
\roman{Hom}(\e\br X/\br C,\Bbb Q/\Bbb Z\e),
$$
where $T$ is the N\'eron-Severi torus of $X$.
\endproclaim

\medpagebreak

When $C=\Bbb P^{1}_{k}$, hypothesis (i) in the statement of
Corollary 2.10 may be replaced by the condition:
``$X\be(k_{v})\neq\emptyset$ for every real prime $v$ of $k$" (see
[6, Corollary 4.6]). Further, $\br X/\e\br
C=H^{@,@,1}\!\left(k,\ns\e\right)$ (see, for example, [3, p.116]).
Thus, we have the following generalization of a well-known theorem
of P.Salberger [10, Theorem 7.5, p.520].

\proclaim{Corollary 2.11} Let $k$ be a number field and let $X\ra
\Bbb P^{1}_{k}$ be a Severi-Brauer $k$-fibration of squarefree
index. Assume that $X\be(k_{v})\neq\emptyset$ for every real prime
$v$ of $k$. Then there exists a natural exact sequence
$$
0\ra{\cyr X}^{\e 1}\be(T\e)\ra A_{0}(X)\ra \bigoplus_{v}
A_{0}(X_{v})\ra\roman{Hom}\!\left(\e
H^{@,@,1}\ng\left(k,\ns\right),\Bbb Q/\Bbb Z\e\right),
$$
where $T$ is the N\'eron-Severi torus of $X$.\qed
\endproclaim

\enddocument

\proclaim{Theorem 8.3} Let $q\:X\ra C$ be a smooth ruled surface
over a number field $k$. Assume: \roster \item"(a)"
$X(k_{v})\neq\emptyset$ for every real prime $v$ of $k$,
\item"(b)" there exists a point $P\in Y(k)$ with smooth fiber
$q^{-1}(P)$, and \item"(c)" the canonical map
$$H^{@,@,1}\be(k,\k2\kb @,@,(X)@,/@,@,H^{@,@,@,0}\be(\xb,\ks2))\ra
\underset{v}\to{\prod}H^{@,@,1}\be(k_{v},\k2\kvb
@,@,(X)@,/@,@,H^{@,@,@,0}\be(\xvb,\ks2))$$ is {\rm{surjective}}.
\endroster
Then there exists a natural exact sequence
%
the map (8.1) (resp. (8.1) for $X=X_{v}$).
\endproclaim
{\it Proof}. See Theorem 6.3 above.\qed

\medpagebreak

{\bf{Remarks 8.4.}} (a) If $q\:X\ra Y$ is {\it any} smooth ruled
surface over a number field $k$, then there exists a canonical map
$\prod_{v}\be\choxvcv\be\ra\ng{\cyr Q}^{@,@,1}\be(S)$ which is
surjective exactly when $\cok\Phi\be\ra\be\prod_{v}\be\cok\Phi_{v}$
is surjective. See Remark 6.4(a).

\smallpagebreak

(b) Regarding hypothesis (c) in the statement of the theorem, see
Remark 6.4(b).

\medpagebreak

We now recall from Section 6 that there exists a natural map
$$\prod_{v}\choxvcv\ra
\roman{Hom}\e(H^{@,@,1}\be(k,\ns),\Bbb Q/\Bbb Z)\tag 8.2$$ which is
surjective exactly when ${\cyr X}^{@,@, 2}(S)=0$ and the natural map
$\cok\Phi\ra\prod_{v}\cok\Phi_{v}$ is surjective (see Remark 6.6(b)
and Examples 4.7 and 5.8).

\medpagebreak

\proclaim{Lemma 8.5} Let $q\:X\ra Y$ be a smooth ruled surface over
a number field $k$. Assume: \roster \item"(a)" there exists a point
$P\in Y(k)$ with smooth fiber $q^{-1}(P)$, and \item"(b)" the
natural map (8.2) is surjective (e.g., $H^{@,@,1}\be(k,\ns)=0$).
\endroster
Then there exists a canonical isomorphism
$${\cyr X}^{@,@,2}(\krn\partial^{\e 1})=
{\cyr X}^{@,@,2}(\k2\kb @,@,(X)@,/@,@,H^{@,@,@,0}\be(\xb,\ks2)).$$
\endproclaim
{\it Proof}. See Lemma 6.7.\qed

\medpagebreak

\proclaim{Lemma 8.6} Let $q\:X\ra Y$ be a smooth ruled surface over
a number field $k$. Assume: \roster \item"(a)" for every real prime
$v$ of $\e k$, $X(k_{v})\neq\emptyset$, and \item"(b)" the natural
map
$$H^{@,@,1}\be(k,\k2\kb @,@,(X)@,/@,@,H^{@,@,@,0}\be(\xb,\ks2))\ra
\underset{v}\to{\prod}H^{@,@,1}\be(k_{v},\k2\kvb
@,@,(X)@,/@,@,H^{@,@,@,0}\be(\xvb,\ks2))$$ is surjective.
\endroster
Then there exists a canonical isomorphism
$${\cyr X}^{@,@,2}(\k2\kb @,@,(X)@,/@,@,H^{@,@,@,0}\be(\xb,\ks2))=
{\cyr X}^{@,@,3}(k(X),\Bbb Q/\Bbb Z(2)).$$
\endproclaim
{\it Proof}. See Lemma 6.8.\qed

\medpagebreak

Combining Lemmas 8.5 and 8.6, we immediately obtain

\proclaim{Proposition 8.7} Let $q\:X\ra Y$ be a smooth ruled surface
over a number field $k$. Assume: \roster \item"(a)"
$X(k_{v})\neq\emptyset$ for every real prime $v$ of $k$,
\item"(b)" there exists a point $P\in Y(k)$ with smooth fiber
$q^{-1}(P)$, \item"(c)" the natural map (8.2) is surjective, and
\item"(d)" the natural map
$$H^{@,@,1}\be(k,\k2\kb @,@,(X)@,/@,@,H^{@,@,@,0}\be(\xb,\ks2))\ra
\underset{v}\to{\prod}H^{@,@,1}\be(k_{v},\k2\kvb
@,@,(X)@,/@,@,H^{@,@,@,0}\be(\xvb,\ks2))$$ is surjective.
\endroster
Then there exists a canonical isomorphism
$${\cyr X}^{@,@,2}(\krn\partial^{\e 1})=
{\cyr X}^{@,@,3}(k(X),\Bbb Q/\Bbb Z(2)).\qed$$
\endproclaim

\medpagebreak

\proclaim{Theorem 8.8} Let $q\:X\ra Y$ be a smooth ruled surface
over a number field $k$. Assume: \roster \item"(a)"
$X(k_{v})\neq\emptyset$ for every real prime $v$ of $k$,
\item"(b)" there exists a point $P\in Y(k)$ with smooth fiber
$q^{-1}(P)$, \item"(c)" the natural map (8.2) is surjective,
\item"(d)" the local periods $\yv$ of $X$ are relatively prime in
pairs, and \item"(e)" the natural map
$$H^{@,@,1}\be(k,\k2\kb @,@,(X)@,/@,@,H^{@,@,@,0}\be(\xb,\ks2))\ra
\underset{v}\to{\prod}H^{@,@,1}\be(k_{v},\k2\kvb
@,@,(X)@,/@,@,H^{@,@,@,0}\be(\xvb,\ks2))$$ is surjective.
\endroster
Then there exists a natural complex
$$\align 0 &\ra\krn\ng\e\Biggl[\cok\pi_{X}^{*}\ra
\prod_{v}\cok\pi_{X_{v}}^{*}\e\Biggr]\ra {\cyr X}^{@,@,3}(k(X),\Bbb
Q/\Bbb Z(2))\ra\\
&\ra{\cyr X}^{@,@,1}\lbe(J_{Y})/N_{2}\,,\endalign$$ where
$\pi_{X}^{*}$ (resp. $\pi_{X_{v}}^{*}$, for each $v$) is the natural
map $\chox\ra\chob^{\g}$ (resp. $\choxv\ra\choxvb^{\gv}$), is the
group (1.2) and $N_{2}$ is a finite group of order
$\y/\e{\tsize{\prod}}\e\yv$. Further, the above complex is exact
except perhaps at , where its homology is a subgroup of
$\e\cok\ng\be\left[\cok\pi_{X}^{*}\ra
\prod_{v}\cok\pi_{X_{v}}^{*}\right]$.
\endproclaim
{\it Proof.} See Corollary 6.10.\qed

\medpagebreak

Recall the group  of rational equivalence classes of 0-cycles of
degree zero on $X$, i.e.,  is the kernel of the degree map
$\:\chox\ra\Bbb Z$. Let
$$\pi_{X,0}^{*}\:\ra J_{Y}\be(k)\tag 8.3$$
be the map which is induced by $\pi_{X}^{*}\:\chox\ra\chob^{\g}$
together with the identification $^{\g}=J_{Y}\be(k)$ [16, p.168].

\proclaim{Corollary 8.9} Let $q\:X\ra Y$ be a smooth ruled surface
over a number field $k$. Assume: \roster \item"(a)" $X$ has a
0-cycle of degree 1, \item"(b)" there exists a point $P\in Y(k)$
with smooth fiber $q^{-1}(P)$, \item"(c)" the natural map (8.2) is
surjective, and \item"(d)" the natural map
$$H^{@,@,1}\be(k,\k2\kb @,@,(X)@,/@,@,H^{@,@,@,0}\be(\xb,\ks2))\ra
\underset{v}\to{\prod}H^{@,@,1}\be(k_{v},\k2\kvb
@,@,(X)@,/@,@,H^{@,@,@,0}\be(\xvb,\ks2))$$ is surjective.
\endroster
Then there exists a natural complex
$$\align 0 &\ra\krn\ng\e\Biggl[\cok\pi_{X,0}^{*}\ra
\prod_{v}\cok\pi_{X_{v}\be,0}^{*}\e\Biggr]\ra {\cyr
X}^{@,@,3}(k(X),\Bbb Q/\Bbb Z(2))\ra\\
&\ra{\cyr X}^{@,@,1}\lbe(J_{Y})\,,\endalign$$ where $\pi_{X,0}^{*}$
(resp. $\pi_{X_{v}\be,0}^{*}$, for each $v$) is the map (8.3) (resp.
(8.3) for $X=X_{v}$) and  is the group (1.2). Further, the above
complex is exact except perhaps at , where its homology is a
subgroup of $\e\cok\ng\be\left[\cok\pi_{X,0}^{*}\ra
\prod_{v}\cok\pi_{X_{v}\be,0}^{*}\right]$.
\endproclaim
{\it Proof}. Hypothesis (a) implies that hypotheses (a) and (d) of
the theorem hold true. Further, the maps $\pi_{X,0}^{*}\:\ra
J_{Y}\be(k)$ and $\pi_{X}^{*}\:\chox\ra\chob^{\g}$ have isomorphic
kernels and cokernels [14, Lemma 4.3]. The corollary is now
immediate from the theorem.\qed

\medpagebreak

{\bf{Questions 8.10.}} Note that $\cok\pi_{X,0}^{*}$ is finitely
generated (by the Mordell-Weil Theorem) and isomorphic to a subgroup
of the torsion group $H^{@,@,2}(k,\krn)$. Hence $\cok\pi_{X,0}^{*}$
is in fact finite (in particular,  is a finitely generated abelian
group of rank equal to the rank of $J_{Y}\be(k)$). Each group
$\cok\pi_{X_{v}\be,0}^{*}$, where $v$ ranges over the set of primes
of $k$, is finite as well. The following questions arise naturally.
\roster \item"Q1." Are the groups $\cok\pi_{X_{v}\be,0}^{*}$ zero
for all but finitely many primes $v$ of $k$?, \item"Q2." What is an
explicit description of the map
$$\krn\ng\e\Biggl[\cok\pi_{X,0}^{*}\ra
\prod_{v}\cok\pi_{X_{v}\be,0}^{*}\e\Biggr]\ra{\cyr
X}^{@,@,3}(k(X),\Bbb Q/\Bbb Z(2))$$ and how can its image be
characterized?
\endroster

\proclaim{Proposition 6.3} Let $k$ be a number field, let $C$ be a
smooth and projective $k$-curve and let $X\ra C$ be a smooth ruled
surface with N\'eron-Severi torus $S$. Assume that
$X_{k_{v}\be(C)}(k_{v}(C))\neq\emptyset$ for all but possibly one
prime $v$ of $k$\footnote{This hypothesis is equivalent to requiring
that $X\ra C$ admit local sections $C_{v}\ra X_{v}$ for all but
possibly one prime $v$ of $k$.}. Then ${\cyr X}^{@,@,2}\be(S\e)=0$.
\endproclaim
{\it Proof}. For each $v$, $X_{k_{v}\be(C)}$ is a curve of genus
$0$, so $\x\!\left(X_{k_{v}\lbe(C)}\right)=1$ or $2$ and
$\x\!\left(X_{k_{v}\lbe(C)}\right)=1$ if and only if
$X_{k_{v}\be(C)}(k_{v}(C))\neq\emptyset$. Thus the hypothesis of the
proposition is equivalent to requiring that the local indices
$\x\!\left(X_{k_{v}\lbe(C)}\right)$ be pairwise coprime. Now Lemma
6.2 shows that ${\cyr
X}^{@,@,1}\be(\pic\e\xb\e/\e\pic^{\circ}\cb\e)=0$. On the other
hand, by [12, Proposition 1.11, p.169$\e$], there is a canonical
isomorphism $\pic^{\circ}\cb=\pic^{\circ}\e\xb$. Now Tate-Nakayama
duality completes the proof.\qed

By [Gros, Prop. 1.11, p.169], there exists a natural exact sequence
$$
0\ra J_{C}\ng\left(\e\kb\e\right)\ra\pic\xb\ra\ns\ra 0.\tag 6
$$

\proclaim{Proposition 1.7} If $k$ is a number field or a local field
of characteristic zero, the exact sequence (6) induces a natural
isomorphism
$$
H^{@,@,1}\ng\left(k,\ns\right)=\krn\ng\left[H^{@,@,2}\be(k,J)\ra
H^{@,@,2}\ng\left(k,\pic\xb\e\right)\right].
$$
\endproclaim
{\it Proof.} We need to show that the map $H^{@,@,1}\be(k,J)\ra
H^{@,@,1}\ng\left(k,\pic\xb\e\right)$ induced by (6) is surjective.
There exists a natural commutative diagram
$$
\CD
\e C @>>>H^{@,@,1}\be(k,J)\\
@VVp^{*}V     @VVV\\
 X @>>>H^{@,@,1}\ng\left(k,\pic\xb\e\right)
\endCD
$$
in which both horizontal maps are surjective (see, e.g., [CT,
p.116]). The left-hand vertical map is surjective as well, by [PS,
Theorem 7.3, p.110], and the proposition follows at once.\qed

\medpagebreak

\proclaim{Corollary 1.5}

\roster

\item"a)" If $v$ is a non-real prime of $k$, then
$H^{@,@,1}\be(k_{v},S)=0$.

\item"b)" If $v$ is a real prime of $k$, then
$H^{@,@,1}\be(k_{v},S)$ is a quotient of $(J(k_{v})/
J(k_{v})^{\circ})^{*}$.

\item"c)" ${\cyr
Q}^{1}\be(S)=H^{@,@,1}\ng\left(k,\ns\e\right)^{*}$.

\endroster
\endproclaim
{\it Proof.} If $v$ is non-real, $H^{@,@,1}\be(k_{v},S)$ is dual to
$H^{@,@,1}\!\be\left(k,\ns\right)$ and $H^{@,@,2}\be(k_{v},J)=0$ by
[Duality, Theorem I.3.2, p.51, and Corollary I.2.4, p.35]. Assertion
(a) now follows at once from the proposition. If $v$ is real, then
$H^{@,@,2}\be(k_{v},J)$ is naturally isomorphic to $J(k_{v})/
J(k_{v})^{\circ}$ and $H^{@,@,1}\be(k_{v},S)$ is again dual to
$H^{@,@,1}\!\be\left(k,\ns\right)$ [op.cit., Remark I.3.7, p.56, and
Theorem I.2.13(b), p.42]. Assertion (b) now follows from the
proposition. By the exactness of (1), assertion (c) is equivalent to
the vanishing of ${\cyr X}^{\e 2}\be(S)$. Now, by [Duality, Theorem
I.6.26(c), p.112], the canonical map
$H^{@,@,2}\be(k,J)\ra\bigoplus_{v\roman{\,real}}
H^{@,@,2}\be(k_{v},J)$ is an isomorphism, whence ${\cyr
X}^{1}\ng\left(\ns\right)=0$ by the proposition. Since the latter
group is dual to ${\cyr X}^{\e 2}\be(S)$ [op.cit., Theorem
I.4.20(a), p.80], the proof is complete.\qed

\enddocument

Let $k$ be a number field and let $C$ be a smooth and projective
$k$-curve. A {\it ruled surface over $C$} is a smooth and projective
$k$-surface $X$ equipped with a dominant $k$-morphism $q\: X\ra C$
such that its generic fiber $X_{k(C)}=X\times_{C}\roman{Spec}\,
k(C)$ is a smooth $k(C)$-curve of genus zero. We will write
$\xb_{\kb(C)}=\xb\times_{\cb}\,\roman{Spec}\, \kb(C)$ and, for each
prime $v$ of $k$, $X_{k_{v}\be(C)}=X_{v}\times_{C_{v}}\roman{Spec}\,
k_{v}(C)$.

Let $q\: X\ra C$ be a ruled surface over $C$ and let
$\overline{q}\:\xb\ra\cb$ be the induced morphism. Then
$\overline{q}$ has only finitely many singular fibers, lying over
points $y_{i}\in C(\kb)$ ($1\leq i\leq r$). Each singular fiber is
the union of two projective lines $L_{i}$ and $\tilde{L}_{i}$ which
meet transversaly and are either fixed or permuted by $\g$. For each
$i$, we will write $y_{i}$ for the class of $y_{i}$ in
$\pic\e\cb=\chocb$. There exists a natural exact sequence of
$\g$-modules
$$
0\ra P_{1}\ra P_{2}\be\oplus\pic\,\cb
\ra\pic\xb\ra\Bbb Z\ra 0\tag 5.1
$$
where $P_{1}=\bigoplus_{i=1}^{r}\be\Bbb Z\e y_{i}$ and
$P_{2}=\bigoplus_{i=1}^{r}\!\big( \Bbb Z\e L_{i}\lbe\oplus\lbe\Bbb
Z\e\tilde{L}_{i}\big)$ are permutation $\g$-modules, the map
$P_{1}\ra P_{2}\oplus\pic\,\cb$ is given by
$y_{i}\mapsto(L_{i}+\tilde{L}_{i},-y_{i})$, and the map
$\pic\xb\ra\Bbb Z$ is given by restricting divisors to the generic
fiber $\xb_{\kb(C)}\simeq\Bbb P^{1}_{\be\kb(C)}$ (cf. [11, proof
of Prop. 1.9, p.169$\e$] and [7, p.201$\e$]). Define the positive
integer $\varepsilon$ by the equality
$$\img\be[\e(\pic\xb)^{\g}\ra\Bbb Z\e]=\varepsilon\e\Bbb Z.\tag 5.2$$
Using the fact that
$\left(\be\pic\xb_{\kb(C)}\be\right)^{\be\g}\!\!\big/\, \pic
X_{k(C)}$ is annihilated by multiplication by the index
$\delta\!\left(X_{k(C)}\right)$ of $X_{k(C)}$, it is not difficult
to see that $\varepsilon$ divides $\delta\!\left(X_{k(C)}\right)$,
which is either 1 or 2 (since $X_{k(C)}$ has genus 0). Now, since
$H^{@,@,1}(k,P_{1})=0$, (5.1) induces a natural exact sequence
$$0\ra P_{1}^{\e\g}\ra P_{2}^{\e\g}\be\oplus(\pic\,\cb\e)^{\g}
\ra(\pic\xb\e)^{\g}\ra\varepsilon\e\Bbb Z\ra 0.\tag 5.3$$ Let $f$ be
the divisor on $\xb$ defined by a smooth fiber of $\overline{q}$.
Then the class of $f$ in $\ns$ (which we also denote by $f$) is
$\g$-invariant (see [7, p.201]). Now, there exists a natural exact
sequence of $\g$-modules
$$0\ra P_{1}\ra P_{2}\be\oplus\Bbb Z f\ra\ns\ra\Bbb Z\ra 0,\tag 5.4$$
where the map $P_{1}\ra P_{2}$ is given by $y_{i}\mapsto
(L_{i}+\tilde{L}_{i},-f)$ and the map $\ns\ra\Bbb Z$ is given by
restricting divisors to the generic fiber (cf. [5, p.106]). Define
the positive integer $\eta$ by the equality
$$\img\be[\e\ns^{\g}\ra\Bbb Z\e]=\eta\e\Bbb Z.\tag 5.5$$
Then $\eta$ divides $\varepsilon$, and there exists a natural exact
sequence
$$0\ra P_{1}^{\e\g}\ra P_{2}^{\e\g}\be\oplus\Bbb Z
f\ra\ns^{\g}\ra\eta\e\Bbb Z\ra 0.\tag 5.6$$ We now combine (5.3) and
(5.6) to obtain the following exact commutative diagram

A two-fold application of the snake lemma to the above diagram
yields the following exact sequence
$$0\ra\Bbb Z/\y\be(C)\e\Bbb Z\ra M'\ra\Bbb Z/(\varepsilon/\eta)\e\Bbb
Z\ra 0,$$ where
$M'=\cok\ng\left[(\pic\xb\e)^{\g}\ra\ns^{\g}\e\right]$ is the group
(3.3) and $\y\be(C)$ is the period of $C$. Similarly, for each prime
$v$ of $k$, there exists a natural exact sequence
$$0\ra\Bbb Z/\yv\be(C)\e\Bbb Z\ra M'_{v}\ra\Bbb
Z/(\varepsilon_{v}/\eta_{\e v})\e\Bbb Z\ra 0.$$ Now consider the
natural exact commutative diagram
$$\minCDarrowwidth{.6cm}
\CD 0 @>>> \Bbb Z/\y(C)\e\Bbb Z @>>> M'
@>>>Z/(\varepsilon/\eta)\e\Bbb
Z @>>>0\\
      @.        @VVV     @VVV     @VVV     @.\\
      0 @>>>\underset{v}\to{\bigoplus}\,\e\Bbb Z/\yv\be(C)\e\Bbb Z
@>>>\underset{v}\to{\bigoplus}\, M'_{v}
@>>>\underset{v}\to{\bigoplus}\,\e\Bbb Z/(\varepsilon_{v}/\eta_{\e
v}\be)\e\Bbb Z @>>> 0\,.\endCD$$ Applying the snake lemma to the
above diagram, we obtain

\proclaim{Lemma 5.1} Let $C$ be a smooth and projective $k$-curve
over a number field $k$ and let $X\ra C$ be a ruled surface over
$C$. Assume: \roster \item"(a)" the local periods $\yv\be(C)$ of $C$
are pairwise coprime, and \item"(b)" there exists a prime $v_{0}$ of
$k$ such that $X_{k_{v}\be(C)}(k_{v}(C))\neq\emptyset$ for all
$v\neq v_{0}$.
\endroster
Then the map (3.6), $M'\ra\oplus M'_{v}$, is surjective. Further,
the integers (3.4) and (3.5), $m'$ and $m'_{v}$, are given by
$$m'=\y\be(C)\be\cdot\be(\varepsilon/\eta)$$
and
$$m'_{v}=\yv\be(C)\be\cdot\be(\varepsilon_{v}/\eta_{v}).$$
where $\varepsilon$ and $\eta$ (resp. $\varepsilon_{v}$ and
$\eta_{\e v}$, for each prime $v$ of $k$) are given by (5.2) and
(5.5) (resp., by $\text{(5.2)}_{v}$ and $\text{(5.5)}_{v}$).
\endproclaim
{\it Proof}. This follows from the preceding discussion and Lemma
1.3 of [9]. Note that hypothesis (b) implies that all except
possibly one of the integers $\varepsilon_{v}/\eta_{v}$ is equal to
one, and therefore that these integers are pairwise coprime.\qed

\medpagebreak

\proclaim{Lemma 5.2} Let $k$ be a number field and let $X\ra C$ be a
ruled surface over a smooth and projective $k$-curve $C$. Assume
that there exists a prime $v_{0}$ of $k$ such that
$X_{k_{v}\be(C)}(k_{v}(C))\neq\emptyset$ for all $v\neq v_{0}$. Then
${\cyr X}^{@,@,1}\be(\ns)=0$.
\endproclaim
{\it Proof}. From the exact sequence (5.4) we obtain a natural exact
sequence
$$0\ra\Bbb Z/\e\eta\e\Bbb Z\ra H^{@,@,1}(k,(P_{2}\oplus\Bbb Z\e
f\e)/P_{\e 1})\ra H^{@,@,1}(k,\ns)\ra 0,$$ where $\eta$ is defined
by (5.5) (as noted previously, $\eta$ divides the index
$\x(X_{k(C)})$). Similar sequences exist over $k_{v}$ for each prime
$v$ of $k$, and we have a natural exact commutative diagram
$$\minCDarrowwidth{.6cm}
\CD 0 @>>> \Bbb Z/\eta\e\Bbb Z @>>>
H^{@,@,1}\!\be\left(k,\frac{P_{2}\oplus\e\Bbb Z f}{P_{\e
1}}\be\right)
 @>>> H^{@,@,1}(k,\ns) @>>>0\\
      @.        @VVV     @VVV     @VVV     @.\\
      0 @>>>\underset{v}\to{\prod}\,\e\Bbb Z/\eta_{\e v}\Bbb Z
@>>>\underset{v}\to{\prod}\,
H^{@,@,1}\!\be\left(k_{v},\frac{P_{2,v}\oplus\e\Bbb Z f_{v}}{P_{\e
1,v}}\be\right)
 @>>>\underset{v}\to{\prod}\,H^{@,@,1}(k_{v},\nsv)@>>> 0\,.\endCD$$
Applying the snake lemma to the above diagram, we obtain a natural
exact sequence
$${\cyr X}^{@,@,1}\be((P_{2}\oplus\e\Bbb Z f\e)/P_{\e 1})
\ra {\cyr X}^{@,@,1}\be(\ns)\ra\cok\!\!\left[\Bbb Z/\eta\e\Bbb
Z\ra\underset{v}\to{\prod}\,\e\Bbb Z/\eta_{\e v}\Bbb Z\right].$$ The
hypothesis implies that the group on the right above vanishes (see
[9, Lemma 1.3]), and it remains only to show that ${\cyr
X}^{@,@,1}\be((P_{2}\oplus\e\Bbb Z f\e)/P_{\e 1})=0$. Since
$H^{@,@,1}(k,P_{2}\oplus\Bbb Z\e f\e)=0$, there exists a natural
injective map
$$H^{@,@,1}(k,(P_{2}\oplus\Bbb Z\e f\e)/P_{\e 1})\hookrightarrow
H^{@,@,2}(k,P_{\e 1})$$ which induces an injection ${\cyr
X}^{@,@,1}\be((P_{2}\oplus\e\Bbb Z f\e)/P_{\e
1})\hookrightarrow{\cyr X}^{@,@,2}\be(P_{\e 1})$. But the latter
group vanishes, by [25, Lemma 1.9, p.19$\e$].\qed

\medpagebreak

We now assume, for simplicity, that the local indices of $X$ are
pairwise coprime. Then so are the local indices and local periods of
$C$. Further, since $\xb=0$ by [11, Prop. 1.12, p.169$\e$], we have
$${\cyr B}_{1}(X)={\cyr B}(X)=\krn\!\!\left[\e
X\ra\bigoplus_{v}\e X_{v}\e\right].$$ Now let $J_{C}$ be the
Jacobian variety of $C$, i.e., $J_{C}=\pic^{\circ}_{C/k}$. Since
$\pic^{\circ}\xb= \pic^{\circ}\cb$ by [11, Prop. 1.11, p.169$\e$],
the following result is an immediate consequence of Theorem 4.4 via
Lemmas 5.1 and 5.2 (see also the proof of Lemma 5.1 and compare
Corollary 4.5).

\proclaim{Theorem 5.3} Let $k$ be a number field, let $C$ be a
smooth and projective $k$-curve and let $X\ra C$ be a ruled surface
over $C$. Assume: \roster \item"(a)" the local indices $\xv\be(X)$
of $X$ are pairwise coprime, \item"(b)" there exists a prime $v_{0}$
of $k$ such that $X_{k_{v}\be(C)}(k_{v}(C))\neq\emptyset$ for all
$v\neq v_{0}$, and \item"(c)" $\e{\cyr X}^{@,@,1}\be(J_{C})$
contains no nonzero infinitely divisible elements.
\endroster
Then there exists a natural exact sequence
$$0\ra{\cyr B}(X)\ra
{\cyr X}^{@,@,1}\be(J_{C}\e)/N_{2}\ra N_{3}\ra 0\,,$$ where $N_{2}$
and $N_{3}$ are groups of orders

where $\varepsilon,\eta_{\e v_{0}},\varepsilon_{v_{0}}$ and $\eta$
are given by (5.2), $\slanted{(5.5)}_{ v_{0}}$,
$\slanted{(5.2)}_{v_{0}}\be$ and (5.5), respectively, $a$ divides
$\prod\!\left(\xv(X)\big/\big[(\pic\xvb\e)^{\gv}\!:\pic\e
X_{v}\big]\right)$ and $\y\be(X)$ is the period of $X$. In
particular, ${\cyr B}(X)$ is finite if and only if $\,{\cyr
X}^{@,@,1}\be(J_{\e C}\e)$ is finite, and in this case
$$ a\e\varepsilon\e\eta_{\e v_{0}}\e\y\be(C)\e\y\be(X)\left[
{\cyr B}(X)\right]=
\varepsilon_{v_{0}}\eta\e(\e\prod\e\xv\be(X),\y\be(X)\e)
\prod\be\yv\be(C)\left[{\cyr X}^{@,@,1}\be(J_{ C})\right].\qed$$
\endproclaim

\medpagebreak

\proclaim{Corollary 5.4} Let $k$ be a number field, let $C$ be a
smooth and projective $k$-curve and let $X\ra C$ be a ruled surface
over $C$. Assume: \roster \item"(a)" $\e{\cyr X}^{@,@,1}\be(J_{
C})=0$, \item"(b)" the local indices $\xv\be(X)$ of $X$ are pairwise
coprime, and \item"(c)" there exists a prime $v_{0}$ of $k$ such
that $X_{k_{v}\be(C)}(k_{v}(C))\neq\emptyset$ for all $v\neq v_{0}$.
\endroster
Then the Hasse principle holds for $\e X$, i.e., the natural map $
X\ra\bigoplus_{v}\! X_{v}$ is injective.
\endproclaim
{\it Proof}. This is immediate from the theorem.\qed

\medpagebreak

\noindent{\bf Remark 5.5.} A formula which is perhaps more appealing
than that of Theorem 5.3 may be obtained as follows. Assume that
${\cyr X}^{@,@,1}\be(J_{C}\e)$ is finite and that conditions (a) and
(b) of Theorem 5.3 hold. Since the local indices of $C$ are pairwise
coprime, Corollary 4.5 applied to $C$ yields an equality
$$\y\be(C)^{2}\e[{\cyr B}(C)]=(\e\prod\e\xv\be(C),\y\be(C)\e)\e
\prod\yv\be(C)\e[{\cyr X}^{@,@,1}\be(J_{C}\e)]$$ (see Remark 2.1 and
Proposition 3.1). Combining the above equality with the formula of
Theorem 5.3, we obtain the identity

where $a'=a\e(\varepsilon\e\eta_{\e
v_{0}}/\varepsilon_{v_{0}}\e\eta)$.

\bigpagebreak

We shall now see that Corollary 3.13 is a particular case of a
significantly more general result (see Corollary 3.17 below). Let
$k$ be a number field and let $X$ be a smooth projective $k$-variety
over a smooth projective $k$-variety $Y$, i.e., there exists a
$k$-morphism $q\:X\ra Y$. Recall the group  of rational equivalence
classes of 0-cycles of degree zero on $\xb$.

\proclaim{Proposition 3.16} Let $k$ be a number field and let
$q\:X\ra Y$ be a smooth projective $k$-variety fibered over a smooth
projective $k$-variety $Y$. Assume: \roster \item"(a)" ${\cyr
X}^{@,@,1}\be(\pic^{\circ}\cb\e)={\cyr X}^{@,@,1}\be(\ns)=0$,
\item"(b)" $\zeta\: M'\ra\bigoplus M'_{v}\e$ is surjective,
\item"(c)" the local indices of $X$ are pairwise coprime, \item"(d)"
the maximal divisible subgroups of ${\cyr X}^{@,@,1}()$ and ${\cyr
X}^{@,@,1}\be(\pic^{\circ}\xb\e)$ are zero. \item"(e)" the canonical
morphism $\overline{q}_{*}\:\ra$ is an isomorphism of $\g$-modules,
and
\item"(f)" for every prime $v$ of $k$, the canonical morphism
$(\overline{q}_{v})_{*}\:\ra$ is an isomorphism of $\gv$-modules.
\endroster
Then the Hasse principle holds for $_{\be\be 1}^{\e\prime}X$.
\endproclaim
{\it Proof}. By (e) and (f), there exists a canonical isomorphism
${\cyr X}^{@,@,1}\be(\e)={\cyr X}^{@,@,1}\be(\e)$. On the other
hand, Roitman's Theorem yields natural isomorphisms ${\cyr
X}^{@,@,1}\be(\aoxb\e)={\cyr X}^{@,@,1}\be(\alb)$ and ${\cyr
X}^{@,@,1}\be(\aoyb\e)={\cyr X}^{@,@,1}\be(\alby)$ (see [12,
Corollary 4.10$\e$]), whence ${\cyr X}^{@,@,1}\be(\alb)={\cyr
X}^{@,@,1}\be(\alby)$. Now there exists a canonical non-degenerate
pairing $({\cyr
X}^{@,@,1}\be(\pic^{\circ}\yb\e)/\e\roman{div})\times({\cyr
X}^{@,@,1}\be(\alby)/\e\roman{div})\ra\Bbb Q/\Bbb Z$ (see the next
section). Thus assumptions (a) and (d) yield ${\cyr
X}^{@,@,1}\be(\alby)=0$, whence ${\cyr X}^{@,@,1}\be(\alb)=0$.
Finally, assumption (d) yields ${\cyr
X}^{@,@,1}\be(\pic^{\circ}\xb\e)=0$, and the proposition follows
from Corollary 3.12.\qed

\medpagebreak

The following result extends Corollary 3.13.

\proclaim{Corollary 3.17} Let $X\ra C$ be an admissible quadric
fibration (in the sense of \,{\rm{[8, Definition 3.1\,]}}), of
relative dimension at least 1, over a smooth projective $k$-curve
$C$ over a number field $k$. Assume: \roster \item"(a)" ${\cyr
X}^{@,@,1}\be(J_{C})={\cyr X}^{@,@,1}\be(\ns)=0$, \item"(b)"
$\zeta\: M'\ra\bigoplus M'_{v}\e$ is surjective, \item"(c)" the
local indices of $X$ are pairwise coprime, and \item"(d)" ${\cyr
X}^{@,@,1}\be(\pic^{\circ}\xb\e)$ contains no nonzero infinitely
divisible elements.
\endroster
Then the Hasse principle holds for $_{\be\be 1}^{\e\prime}X$.
\endproclaim
{\it Proof}. Indeed, $X$ satisfies assumptions (e) and (f) of the
proposition. See [8, Th\'eor\`eme 6.1 and proof of Th\'eor\`eme
6.2\,].\qed

\medpagebreak

 The following key result is due to
J.-L.Colliot-Th\'el\`ene [4].

\proclaim{Theorem 1.2} Let $\alpha\in H^{@,@,1}\be(k, S)$. If, for
each prime $v$ of $k$, there exists $z(v)\in\choxvcv$ such that
$\alpha_{v}=\Phi_{v}(z(v))$, then there exists $z\in\choxc$ such
that $\Phi(z)=\alpha$ and $z_{v}=z(v)$ for each $v$.
\endproclaim
{\it Proof}. See [4, Th\'eor\`eme 1.3].\qed

\medpagebreak

\proclaim{Corollary 1.3} Let $\Phi$ be the map (4) and, for each
prime $v$ of $\e k$, let $\Phi_{v}$ the corresponding map associated
to the $k_{v}$-variety $X_{v}$. Then

\roster

\item"a)" The natural map $\krn\Phi\ra\prod_{v}\krn\Phi_{v}$ is
surjective.

\item"b)" The natural map $\cok\Phi\ra\prod_{v}\cok\Phi_{v}$ is
injective.

\endroster
\endproclaim
{\it Proof}. For (a), let $(z(v))\in\prod_{v}\krn\Phi_{v}$ and take
$\alpha=0$ in the theorem. Statement (b) is immediate from the
theorem.\qed

\enddocument